\documentclass{amsart}
\usepackage{amssymb, amsmath, amsfonts, amscd}

\input xypic

\theoremstyle{plain}

\newtheorem{theorem}{Theorem}[section]

\newtheorem{lemma}[theorem]{Lemma}
\newtheorem{proposition}[theorem]{Proposition}

\newtheorem{example}[theorem]{Example}

\theoremstyle{definition}
\newtheorem{definition}[theorem]{Definition}

\theoremstyle{remark}

\numberwithin{equation}{theorem}

\renewcommand{\O}{\mathcal{O} }
\renewcommand{\a}{\alpha}

\renewcommand{\P}{\mathbf{P} }

\newcommand{\pa}{\partial}

\newcommand{\Der}{\operatorname{Der} }

\newcommand{\End}{\operatorname{End} }

\renewcommand{\H}{\operatorname{H} }

\newcommand{\D}{\operatorname{D} }

\newcommand{\C}{\operatorname{C} }

\newcommand{\K}{\operatorname{K_0} }
\newcommand{\hK}{\operatorname{K} }
\newcommand{\Kl}{\operatorname{K_0^\lambda}}

\renewcommand{\lg}{\mathfrak{g}}
\newcommand{\lh}{\mathfrak{h}}
\renewcommand{\ll}{\mathfrak{l}}

\newcommand{\sym}{\operatorname{Sym}}

\newcommand{\w}{\omega}

\newcommand{\Z}{\mathbb{Z}}

\newcommand{\kp}{K^{1/p}}
\def\mod{\text{{\bf conn}-}}
\def\modpf{\text{{\bf mod}-}}
\def\modf{\text{{\bf flat}-}}

\newcommand{\lring}{\underline{$\lambda$-Rings}}
\newcommand{\ring}{\underline{Rings}}
\newcommand{\Vect}{\operatorname{Vect}}

\renewcommand{\C}{\underline{C}}
\renewcommand{\D}{\underline{D}}
\newcommand{\one}{\mathbf{1}^\otimes}
\newcommand{\cat}{\underline{Cat}^{\otimes}_{\wedge}}
\begin{document}

\title{On operations and characteristic classes} 
\author{Helge Maakestad }
\address{NTNU, Trondheim} 
\email{Helge.Maakestad@math.ntnu.no}

\date{March 2009} 

\begin{abstract} We define for any abelian category with
 tensor and exterior products the grothendieck ring
 functor. Furthermore we use
 exterior products to define gamma operations and virtual Chern classes
and virtual Segre classes for arbitrary elements in the grothendieck ring. We apply this to
  define characteristic classes with values in algebraic K-theory and and K-theory of connections.
\end{abstract}

\maketitle
\tableofcontents

\section{Introduction}

The aim of this paper is to give a direct and elementary construction of
virtual characteristic classes with values in the K-theory of an
arbitrary abelian category with tensor and exterior products using operations.
We apply this to construct characteristic classes
with values in algebraic K-theory and K-theory of connections.

To construct characteristic classes of locally free sheaves with
values in algebraic and topological $\hK$-theory one use the \emph{projective bundle
  formula}. In the more general case of the $\hK$-theory of an abelian category with
tensor and exterior product there is no such formula available. To
define characteristic classes using exterior products is an
alternative approach to the construction of characteristic classes
not relying on the projective bundle formula.

The main results of the paper are the following theorems:
Let $\cat$ be the category of abelian categories with tensor and
exterior products and morphisms. Let \lring \text{ } denote the
category of finite dimensional augmented $\lambda$-rings and morphisms.

\begin{theorem}There is a covariant functor
\[ \Kl:\cat \rightarrow \text{\lring} \]
defined by
\[ \Kl(\C)=\{ \K(\C),\{ \lambda^n \}_{n\geq 0}\} .\]
Here $\K(\C)$ is the grothendieck ring of $\C$ and $\lambda^n$ is the
$n$'th exterior product.
\end{theorem}

Let $\C\in Ob(\cat)$. We define for any $l\geq 0$ and $x\in \K(\C)$ the \emph{virtual
characteristic class} $c_l(x)\in \K(\C)$.
\begin{theorem} Let $x,y\in \K(\C)$. The following
  formulas hold:
\begin{align}
&c_l(x+y)=\label{cclass1}\sum_{i+j=l}c_i(x)c_j(y) \\
&c_l(x)=0\text{ if $x$ is effective and $l >e(x)$}\\
&f^*(c_l(x))=c_l(f^*(x))
\end{align}
Here $f:\D\rightarrow \C$ is a morphism in $\cat$.
\end{theorem}

\section{Exterior products and characteristic classes}

In this section we prove existence of virtual Chern and Segre classes
$c_i(x),s_i(x)\in \K(\C)$ for any integer $i\geq 0$ and any $x\in \K(\C)$
where $\C$ is an arbitrary abelian category with tensor and exterior
products.

Let $\C$ be a small abelian category. We say that $\C$ is an
\emph{abelian category with tensor and exterior products} (for short
ACTEP) if the following
holds.
There is a \emph{tensor product} 

\[ \otimes : \C \times \C \rightarrow \C \]

satisfying the following properties:
The triple $\{\C,\oplus,\otimes\}$ is an abelian tensor category.
There are  canonical isomorphisms
\[ x\otimes y \cong y\otimes x \]
and
\[ (x\oplus y)\otimes z\cong x\otimes z\oplus y\otimes z\]
for any objects $x,y,z\in Ob(\C)$. There is a unique object $\one \in
Ob(\C)$ and canonical isomorphisms
\[ \one \otimes x \cong x\otimes \one \cong x \]
for any object $x\in Ob(\C)$. Moreover the endofunctor $x\otimes -$ is
exact for any object $x\in Ob(\C)$.

There is a function
\[ rk:Ob(\C)\rightarrow \mathbb{N} \]
such that $rk(\one)=1$ and $rk$ is additive on exact sequences. 
There is for every $k\geq 0$ an endofunctor - \emph{exterior product} 
\[ \lambda^k:\C \rightarrow \C \]
with the following properties: $\lambda^i(x)=0$ for $i>rk(x)$,
$\lambda^0(x)=x$ for all $x\in Ob(\C)$.
Moreover, for any exact sequence 
\[ 0 \rightarrow x'\rightarrow x \rightarrow x'' \rightarrow 0\]
in $\C$ and any $p\geq 2$ there is a filtration
\[ 0=F_{p+1}\subseteq F_p\subseteq \cdots \subseteq F_1\subseteq
F_0=\lambda^p(x) \]
and canonical isomorphisms
\[ F_i/F_{i+1}\cong \lambda^i(x')\otimes \lambda^{p-i}(x'') \]
in $\C$.
We say a functor
\[ f:\D\rightarrow \C \]
between ACTEP's is a
\emph{morphism of ACTEP's} if
$f$ is an additive tensor functor commuting with the exterior
product. This means for any $l\geq 0$ there are canonical isomorphisms
\[ f(\lambda^l(x))\cong\lambda^l(f(x)) \]
in $Ob(\C)$. 

\begin{definition}
Let $\cat$ denote the category with  ACTEP's 
as objects and morphisms of ACTEP's as morphisms.
\end{definition}

\begin{theorem} There is a covariant functor
\[ \K:\cat \rightarrow \ring \]
defined by
\[ \C \rightarrow \K(\C) \]
where $\K(\C)$ si the grothendieck ring of the category $\C$.
\end{theorem}
\begin{proof} Let for any $\C\in Ob(\cat)$ $\K(\C)$ be the
  grothendieck ring of the category $\C$. Direct sum induce an
  addition operation 
  and tensor product induce a multiplication with $[\one]$ as
  multiplicative unit. One checks that for any morphism
  $f:\D\rightarrow \C$ there is an induced map of rings
\[ f^*:\K(\D)\rightarrow \K(\C) .\]
Finally for two composable morphisms $f,g$ it follows $(g\circ
f)^*=g^*\circ f^*$ and the Theorem follows.
\end{proof}

Define for any element $x\in Ob(\C)$ and any integer $l\geq 0$
$\lambda^l[x]=[\lambda^l(x)]$ where $\lambda^l$ is the $l'th$ exterior
power. It follows that for any exact sequence
\[ 0\rightarrow x' \rightarrow x\rightarrow x'' \rightarrow 0\]
in $\C$ there is for all $p\geq 2$ an equality
\[ \lambda^p[x]=\sum_{i+j=p}\lambda^i[x']\lambda^j[x'']\]
in $\K(\C)$.
Let $\K(\C)[[t]]$ be the ring of formal powerseries in $t$ with
coefficients in $\K(\C)$. Let $1+\K(\C)[[t]]$ be the multiplicative
subgroup of $\K(\C)[[t]]$ consisting of formal powerseries with
constant term equal to one. Let $\Phi(\C)$ denote the abelian monoid
on $\C$ with direct sum as addition operation. Define the following
map
\[ \lambda_t:\Phi(\C)\rightarrow 1+\K(\C)[[t]] \]
by
\[ \lambda_t(x)=\sum_{l\geq 0}\lambda^l(x)t^l.\]
For any exact sequence 
\[ 0\rightarrow x' \rightarrow x \rightarrow x'' \rightarrow 0\]
in $\C$ we get an equality of formal powerseries
\[ \lambda_t(x)=\lambda_t(x')\lambda_t(x'').\]
We get a well defined map of abelian groups
\[ \lambda_t:\K(\C)\rightarrow 1+\K(\C)[[t]] \]
defined by
\[ \lambda_t(n[x]-m[y])=\lambda_t(x)^n\lambda_t(y)^{-m}.\]
We can define for any element $\omega\in \K(\C)$ a formal powerseries
\[ \lambda_t(\omega)=\sum_{l\geq 0}\lambda^l(\omega)t^l .\]
Let $Op(\K)$ be the set of natural transformations of the underlying
set valued functor of $\K$. It follows $\lambda^n\in Op(\K)$

\begin{lemma}The set $Op(\K)$ is an associative ring.
\end{lemma}
\begin{proof}The proof is left to the reader as an exercise.
\end{proof}
Let \lring  \text{ }denote the category of finite dimensional augmented
$\lambda$-rings.

\begin{theorem}There is a covariant functor
\[ \Kl:\cat \rightarrow \text{\lring} \]
defined by 
\[ \Kl(\C)=\{ \K(\C),\{ \lambda^n \}_{n\geq 0}\} .\]
Here $\K(\C)$ is the grothendieck ring of $\C$ and $\lambda^n$ is the
operation defined above.
\end{theorem}
\begin{proof} One checks that for any category $\C$ $\{\K(\C),\{
  \lambda^n\}_{n\geq 0} \}$ is a $\lambda$-ring. Any morphism
  $f:\D\rightarrow \C$ induce a morphism
  $f^*:\K^{\lambda}(\D)\rightarrow \K^{\lambda}(\C)$ of
  $\lambda$-rings because $f$ commutes with exterior
  products. Moreover since any object in $\C$ has finite rank it
  follows that $\K^{\lambda}(\C)$ is a finite dimensional
  $\lambda$-ring. Finally the rank function $rk$ defines a map
\[ rk :\K(\C)\rightarrow \mathbb{Z} \]
which is a map of $\lambda$-rings where $\mathbb{Z}$ has the canonical
$\lambda$-structure, and the Theorem is proved.
\end{proof}

Let $u=t/1-t$ and let $\gamma_t(x)=\lambda_u(x)=\sum_{l\geq
  0}\gamma^l(x)t^l$.
It follows that for any exact sequence in $\C$
\[ 0\rightarrow x' \rightarrow x \rightarrow x'' \rightarrow 0\]
there is an equality
\[ \gamma_t(x)=\gamma_t(x')\gamma_t(x'').\]
It follows that for any $l\geq 0$ we have $\gamma^l\in Op(\K)$.
\begin{definition} Let $\omega=\sum_i n_i[x_i]\in \K(\C)$. Define
\[ e:\K(\C)\rightarrow \mathbb{Z} \]
by
\[ e(\omega)=\sum_{i} n_i rk(x_i).\]
Define
\[ d:\K(\C)\rightarrow \K(\C) \]
by
\[ d(\omega)=e(\omega)[\one] \]
where $\one$ is the unit object for $\otimes$.
Define for all $l\geq 0$
\[ c_l(\omega)=(-1)^l\gamma^l(\omega-d(\omega)) \]
to be the $l$'th characteristic class of $\omega\in \K(\C)$.
We say an element $\omega$ is effective if $n_i\geq 1$ for all $i$.
\end{definition}

\begin{theorem} \label{maintheorem} Let $x,y\in \K(\C)$. The following
  formulas hold:
\begin{align}
&c_l(x+y)=\label{cclass1}\sum_{i+j=l}c_i(x)c_j(y) \\
&c_l(x)=\label{cclass2}0\text{ if $x$ is effective and $l >e(x)$}\\
&f^*(c_l(x))=\label{cclass3}c_l(f^*(x))
\end{align}
Here $f:\D\rightarrow \C$ is a morphism in $\cat$.
\end{theorem}

\begin{proof}

\end{proof}
\begin{proof}We prove \ref{cclass1}:
\[
c_l(x+y)=(-1)^l\gamma^l(x+y-d(x+y))=(-1)^l\gamma^l(x-d(x)+y-d(y))=\]
\[(-1)^l\sum_{i+j=l}\gamma^i(x-d(x))\gamma^j(y-d(y))=\sum_{i+j=l}(-1)^i\gamma^i(x-d(x))(-1)^j\gamma^j(y-d(y))=\]
\[\sum_{i+j=l}c_i(x)c_j(y),\]
and \ref{cclass1} is proved.

We next prove \ref{cclass2}: Note: $u=t/1-t$. We see that
\[ \gamma_t(\one)=\lambda_u(\one)=1+[\one]u=1+t/1-t=1/1-t.\]
We get:
\[
\gamma_t(x-d(x))=\gamma_t(x)\gamma_t(\one)^{-e(x)}=\gamma_t(x)(1-t)^{e}\]
where $e=e(x)$ and $x=\sum_i n_i[x_i]=[\oplus_i x_i^{n_i}]=[\omega] \in
\K(\C)$.
We get:
\[\gamma_t(x-d(x))=\sum_{l\geq
  0}\gamma^l(x-d(x))=\gamma_t(x)(1-t)^e=\]
\[\lambda_u(\omega)(1-t)^e=(1+\lambda^1(\omega)u+\lambda^2(\omega)u^2+\cdots
+\lambda^e(\omega)u^e)(1-t)^e=\]
\[(1-t)^e+\lambda^1(\omega)t(1-t)^{e-1}+\cdots +\lambda^e(\omega)t^e.\]
It follows that
\[\gamma_t(x-d(x))= \sum_{l\geq 0}\gamma^l(x-d(x))=p_0+p_1t+\cdots
+p_et^e\]
hence $\gamma^l(x-d(x))=0$ for $l\geq e=e(x)$ and hence
\[ c_l(x)=(-1)^l\gamma^l(x-d(x))=0 \]
for $l> e(x)$, and \ref{cclass2} is proved.

We next prove \ref{cclass3}:
\[
f^*(c_l(x))=f^*((-1)^l\gamma^l(x-d(x)))=(-1)^l\gamma^l(f^*(x)-f^*(e(x)[\one]))=\]
\[ (-1)^l\gamma^l(f^*(x)-e(f^*(x))[\one])=c_l(f^*(x)),\]
and \ref{cclass3} is proved.
\end{proof}

\begin{definition}Let $x\in \K(\C)$ be an element. Let
\[ c_t(x)=\sum_{k\geq 0}c_k(x)t^k\in 1+\K(\C)[[t]] \]
be the \emph{Chern power series} of the element $x$.
\end{definition}

It follows from Theorem \ref{maintheorem} that for any functor
$f:\D\rightarrow \C$ in $\cat$ there is a commutative diagram
\[
\diagram  \K(\D) \dto^{f^*} \rto^{c_t} & 1+\K(\D)[[t]] \dto^{1\otimes
  f^*} \\
          \K(\C) \rto^{c_t} & 1+\K(\C)[[t]].
\enddiagram \]
Hence the Chern power series define a natural transformation
\[ c_t:\K(-) \rightarrow  1+\K(-)[[t]] \]
of functors.
Here we view $\K(-)$ and $1+\K(-)[[t]]$ as functors
\[ \K(-),1+\K(-)[[t]]:\cat \rightarrow \underline{Ab} \]
where $\underline{Ab}$ is the category of abelian groups.

\begin{definition} We say that a natural transformation
\[ u: \K(-)\rightarrow 1+\K(-)[[t]] \]
of functors is a \emph{theory of characteristic classes with values in $\K$}.
Let $Class(\K)$ denote the \emph{set of theories of  characteristic classes
with values in $\K$}.
\end{definition} 

It follows the Chern power series define a theory of characteristic
classes $c_t\in Class(\K)$.

Since the Chern power series $c_t(x)$ is a unit for
the multiplication in $\K(\C)[[t]]$ there exists a power series 
$s_t(x)\in 1+\K(\C)[[t]]$ - the \emph{Segre
  power series} - defined by  $s_t(x)=c_t(x)^{-1}$. Let 
\[ s_t(x)=\sum_{k\geq 0}s_k(x)t^k .\]
We get
\[ s_t(x+y)=c_t(x+y)^{-1}=(c_t(x)c_t(y))^{-1}=c_t(x)^{-1}c_t(y)^{-1}=s_t(x)s_t(y)\]
and
\[f^*(s_t(x))=f^*(c_t(x)^{-1})=(f^*(c_t(x)))^{-1}=c_t(f^*(x))^{-1}=s_t(f^*(x)).\]
\begin{definition} We let the classes $s_k(x)\in \K(\C)$ for $k\geq 0$
  be the \emph{virtual Segre classes} of the element $x\in \K(\C)$.
\end{definition}
It follows the Segre classes satisfy the following formulas:
Let $x,y\in \K(\C)$ be arbitrary elements and let $f:\D\rightarrow \C$
be a morphism in $\cat$.

\begin{align} 
s_k(x+y)=& \label{segre1}\sum_{i+j=k}s_i(x)s_j(y) \\
f^*(s_k(x))=& \label{segre2} s_k(f^*(x))
\end{align}
Hence the Segre power series define a natural transformation
\[ s_t(-):\K(-)\rightarrow 1+\K(-)[[t]] \]
of functors. It follows $s_t\in Class(\K)$
It is not true in general that $s_i(x)=0$ for $x$ effective and $i>e(x)$.

\begin{example} Characteristic classes in algebraic and topological $\hK$-theory.
\end{example}
We include a discussion of existence of characteristic classes with
values in algebraic and topological grothendieck groups using the
constructions above.

Following \cite{sga6}, let $X$ be an arbitrary scheme and let $vb(X)$ be the category of
locally free finite rank $\O_X$-modules. It follows that $vb(X)$ is an ACTEP.

It is a standard
fact that the grothendieck ring $\K(X)=\K(vb(X))$ of $vb(X)$ is a commutative
ring with unit. Direct sum induce the addition operation and tensor
product induce the  multiplication. 

Given any morphism $f:Y\rightarrow X$ the pull-back $f^*$ defines a
map of rings
\[ f^*:\K(X)\rightarrow \K(Y) .\]
Let 
\[ 0 \rightarrow E' \rightarrow E \rightarrow E'' \rightarrow 0 \]
be an exact sequence of locally free sheaves of ranks $e', e$ and $e''$.
From \cite{hartshorne}, Exercise II.5.16,  there is for all $2\leq l \leq e$ a filtration
\begin{align} \label{filtr}
 0&=F_{l+1}\subseteq F_l \subseteq \cdots \subseteq F_0=\wedge^l E 
\end{align}
where 
\[ F_i/F_{i+1}\cong \wedge^i(E')\otimes \wedge^{l-i}(E'') \]
for all $0 \leq i \leq l$.
Let $Op(\K)$ be the set of natural transformations of the underlying
set-valued functor of the functor
$\K$. The operation $\lambda^k[E]=[\wedge^k E]$
extends to give a natural transformation $\lambda^k \in Op(\K)$ for
all $k\geq 0$. Let $\underline{Sch}$ denote the category of schemes. The construction
$\K$ and $\lambda^n$ defines a contravariant functor
\[ \Kl:\underline{Sch}\rightarrow \text{\lring} \]
by
\[ \Kl(X)=\{\K(vb(X)), \{\lambda^n\}_{n\geq 0} \} .\]

Let $u=t/1-t$ and define the powerseries
\[ \gamma_t(x)=\lambda_u(x)=\sum_{l\geq 0}\lambda^l(x)u^l=\sum_{l\geq
  0}\gamma^l(x)t^l\in 1+\K(X)[[t]].\]
Following \cite{sga6} (and \cite{karoubi}) 
we get operations - \emph{gamma operations} - $\gamma^l\in Op(\K)$ for all $l\geq 0$. 

\begin{definition} \label{class2} Let $x=\sum_in_i[E_i]\in
  \K(X)$ be an element.
Define the following map
\[ e:\K(X)\rightarrow \mathbb{Z} \]
by
\[ e(x)=\sum_i n_i rk(E_i).\]
Define furthermore
\[ d:\K(X)\rightarrow \K(X) \]
by
\[ d(x)=e(x)[\O_X] .\]
For any integer $l\geq 0$ define the \emph{l'th characteristic class}
of $x$ to be
\[ c_l(x)=(-1)^l\gamma^l(x-d(x)) \in \K (X).\]
\end{definition}

\begin{proposition}\label{properties2} Let $f:Y\rightarrow X$ be a map of
  schemes and let $x,y\in \K(X)$ be arbitrary elements.
The following formulas hold:
\begin{align}
c_l(x+y)=\label{2cclass1} &\sum_{i+j=l}c_i(x)c_j(y) \\
c_l(x)=&\label{2cclass2} 0 \text{ if }l> e(x)\\
f^*(c_l(x))=\label{2cclass3}&c_l(f^*(x))
\end{align}
\end{proposition}
\begin{proof} The proof  is similar to the proof of Theorem
  \ref{maintheorem} and is left to the reader as an exercise.
\end{proof}
Let $c_t(x)=\sum{k\geq 0}c_k(x)t^k\in 1+\K(X)[[t]]$ be the Chern power
series of $x\in \K(X)$. Define $s_t(x)=c_t(x)^{-1}$. It follows we get
Segre classes $s_k(x)\in \K(X)$ with $k\geq 0$ for any element $x=\sum_i n_i[E_i]$.

If one considers the grothendieck ring $\K(B)$ of finite rank
continuous vector bundles on a topological space $B$ a similar
construction using tensor operations defines characteristic classes
$c_i(E)\in \K(B)$ where $E$ is a complex continuous vector bundle on
a topological space $B$.
If $E$ is a rank $n$ vector bundle on $B$ it follows we get
a \emph{total characteristic class} 
\[ c(E)=\sum_{i= 0}^n c_i(E)\in \K(B) .\]
Let 
\[ \pi: \P(E^*)\rightarrow B \]
be the associated projective bundle of $E$. Its fiber over a point
$x\in B$ is the projectivization $\P(E_x^*)$ of the dual vector space
$E_x^*$ where $E_x$ is the fiber of $E$ at $x$. By \cite{karoubi} Proposition IV.7.4 the following holds:
The map
\[ \pi^*: \K(B)\rightarrow \K(\P(E^*)) \]
is injective and the ring $\K(\P(E^*))$ is free of rank $n$ on the
element
\[ h=1-[\O(-1)] \]
where $\O(-1)$ is the tautological bundle on $\P(E^*)$. It follows we
get an equation
\begin{align}\label{equation}
& h^n-c_1(E)h^{n-2}+c_2(E)h^{n-2}+\cdots +(-1)^nc_n(E)=0 
\end{align}

where the classes $c_i(E)$ in Equation \ref{equation} are the ones
defined by  operations on $\K$. 
The classes defined by Equation \ref{equation} are the \emph{Chern-classes} of the
vector bundle $E$ with values in $\K(B)$. Defining the Chern-classes
of a vector bundle using the projective bundle $\P(E^*)$ is usually
referred to as the \emph{projective bundle formula}. 

The characteristic class $c_i(E)\in \K(B)$ is related to the
exterior product of bundles in the following way (see \cite{karoubi} Section IV.2.18)
\[ c_i(E)=\binom{n}{i}[\wedge^0 E]-\binom{n-1}{i-1}[\wedge^1 E]+\cdots
+(-1)^i\binom{n-i}{0}[\wedge^i E],\]
hence we may use the exterior product to define well behaved
characteristic classes with values in  $\K$-theory. 

Defining the
characteristic classes $c(x)$ directly using operations
is an alternative approach to the theory of
characteristic classes going ``around'' the projective bundle formula.
As we have seen: In the case of the category $\cat$
there is no replacement for the projective bundle
formula. Hence operations give a direct, intrinsic and
elementary approach to the construction of characteristic classes.

\section{Characteristic classes of connections}

In this section we define characteristic classes with values in $\K(\mod \lg )$ where $\lg$ is
a restricted Lie-Rinehart algebra using exterior products. We introduce
$\lambda$-operations on the grothendieck ring $\K(\mod \lg )$ using
techniques similar to the ones found in \cite{karoubi} and \cite{kratzer}.

Let $k\subseteq K$ be fields of characteristic $p>0$. 
\begin{definition} A sub $k$-Lie-algebra and $K$-vector space $\lg \subseteq \Der_k(K)$ is a 
\emph{$p$-$(k,K)$-Lie algebra } if for any $\pa \in \lg$ it follows
that $\pa^{[p]}=\pa \circ \cdots \circ \pa \in \lg$. Given two $p$-$(k,K)$-Lie algebras $\lg$ and $\lh$ a
\emph{morphism} of $p$-$(k,K)$-Lie algebras is given by an an inclusion  $I:\lg \rightarrow
\lh$. Let $Lie_{K/k}$ denote the category of $p$-$(k,K)$-Lie algebras and morphisms.
A \emph{connection} $\rho$ is a map of $K$-vector spaces
\[ \rho:\lg\rightarrow \End_k(V) \]
where $V$ is a finite dimensional $K$-vector space satisfying the
following formula:
\[ \rho(\pa)(ax)=a\rho(\pa)(x)+\pa(a)x  \]
for all $a\in K$ and $x\in V$. 
The \emph{curvature} $R_\rho$ of $\rho$ is the map 
\[ R_\rho(\pa, \pa')= \rho([\pa, \pa'])-[\rho(\pa),\rho(\pa')] .\]
The map 
\[ \psi_\rho(\pa)=\rho(\pa^{[p]})-\rho(\pa)^{[p]} \]
is the \emph{p-curvature} of $\rho$.
The connection is \emph{flat} if $R_\rho=0$. It is
\emph{p-flat} if $R_\rho=\psi_\rho=0$. 
The connection $(V,\rho)$ is \emph{nilpotent of exponent $\leq n$} if there
is a filtration of $\lg$-connections
\[ 0=F^n\subseteq F^{n-1}\subseteq \cdots \subseteq F^1 \subseteq
F^0=V \]
where the induced connection $F^i/F^{i+1}$ has $p$-curvature zero for
all $i$. 
Let $(V,\rho)$ and $(W,\eta)$ be $\lg$-connections. A $K$-linear
  map
\[ \phi: V\rightarrow W \]
is a map of $\lg$-connections if for all $x\in \lg$ there is a
commutative diagram
\[ \diagram   V\rto^{\phi}\dto^{\rho(x)} & W \dto^{\eta(x)} \\
              V \rto^{\phi} & W .
\enddiagram \]
Let $\mod\lg $ (resp. $\modf \lg$, $\modpf \lg$) denote the category
of $\lg$-connections (resp flat $\lg$-connections, p-flat
$\lg$-connections) of finite
dimension over $K$ and morphisms
of $\lg$-connections. 
\end{definition}

\begin{theorem} There is a contravariant functor
\[ \K: Lie_{K/k}\rightarrow \ring\]
where for all $\lg\in Lie_{K/k}$ $\K(\lg)$ is the grothendieck ring of
the category $\mod \lg)$.
\end{theorem}
\begin{proof}The proof is left to the reader as an exercise.
\end{proof}

Note: A $p$-$(k,K)$-Lie algebra is also referred to as a
\emph{restricted Lie-Rinehart algebra}. A $p$-flat $\lg$-connection is
also referred to as a \emph{restricted $\lg$-module}. A connection
$\rho:\lg \rightarrow \End_k(V)$ is flat if and only if $\rho$ is a
morphism of $k$-Lie algebras.
There exists a restricted
enveloping algebra $U^{[p]}(\lg)$ of the restricted Lie-Rinehart algebra $\lg$ with
the following property: There is an equivalence of categories
\[ \modpf \lg \cong \modpf U^{[p]}(\lg)  ,\]
where $\modpf U^{[p]}(\lg)$ is the category of finite dimensional left $U^{[p]}(\lg)$-modules.

\begin{lemma}Let $(V,\rho)$ be a $\lg $-connection, where 
 $\lg \subseteq \Der_k(K)$ is a restricted Lie-Rinehart algebra.
It follows $V^{\otimes l},\wedge^l V$ and $\sym^l(V)$ are $\lg$-connections.
\end{lemma}
\begin{proof}The proof is left to the reader as an exercise.
\end{proof}

Let
\[ 0\rightarrow U \rightarrow V \rightarrow W \rightarrow 0 \]
be an exact sequence in $\mod \lg$, where $\lg \subseteq
\Der_k(K)$ is a restricted Lie-Rinehart algebra.

There is the following Lemma:
\begin{lemma} \label{filtration} There is for all $l\geq 2$ a
  $\lg$-stable filtration
\[0= F_{l+1}\subseteq F_l \subseteq \cdots \subseteq F_1 \subseteq
F_0=\wedge^l V \]
where $F_i$ is a $\lg$-connection, with the following property:
\[  F_i/F_{i+1}\cong \wedge^i U\otimes \wedge^{(l-i)}W \]
is an isomorphism of $\lg$-connections.
\end{lemma}

\begin{proof} Let $\pi :V^{\otimes l}  \rightarrow \wedge^l V$ be the
  canonical projection map. It follows $\pi$ is a map of
  $\lg$-connections. There is a filtration of $\lg$-connections
\[ 0\subseteq U^{\otimes l}\subseteq U^{\otimes (l-1)}\otimes V \subseteq \cdots
\subseteq U^{\otimes (l-i)}\otimes V^{\otimes i} \subseteq \cdots
\subseteq V^{\otimes l}.\]
Make the following definition:
\[ F_i=\pi (U^{\otimes i}\otimes V^{\otimes (l-i)} )\subseteq \wedge^l
V.\]
Since $\pi$ is a map of $\lg$-connections, it follows that the filtration
\[0= F_{l+1}\subseteq F_l \subseteq \cdots \subseteq F_1 \subseteq
F_0=\wedge^l V \]
is a filtration of $\lg$-connections.
There is a commutative diagram of $\lg$-connections
\[ 
\diagram  U^{\otimes i}\otimes V^{\otimes (l-i)}\dto^{1\otimes p}
\rto^{\pi} & F_i \dto \\
          U^{\otimes i}\otimes W^{\otimes (l-i)} \rto^{\tilde{\pi}} \dto &
          F_i/F_{i+1}\dto  \\
         \wedge^i U\otimes \wedge^{(l-i)} W \rto^g & F_i/F_{i+1} 
\enddiagram.
\]
The claim is that the bottom horizontal map
\[ g:\wedge^i U\otimes \wedge^{(l-i)} W \rightarrow F_i/F_{i+1}  \]
is an isomorphism of $\lg$-connections. It is enough to prove it is an
isomorphism of $K$-vector spaces.

The sequence
\[ 0\rightarrow U \rightarrow V \rightarrow W \rightarrow 0 \]
is split as sequence of $K$-vector spaces, hence $V=U\oplus W$ as
$K$-vector space. Consider the diagram
\[
\diagram   U^{\otimes i}\otimes V^{\otimes (l-i)}\rto \dto &
V^{\otimes l} \dto \\
              F_i \rto & \wedge^l V   
\enddiagram.
\]
Since $V=U\oplus W$ as $K$-vector space the following holds:
\[ F_i =f(U^{\otimes i}\otimes V^{\otimes (l-i)})=\wedge^i U\wedge
(\wedge^{(l-i)}(U\oplus W))=\]
\[\wedge^i U\wedge \oplus_{j=0}^{l-i}\wedge^j U\otimes \wedge^{l-i-j}W
=\]
\[ \oplus_{j=0}^{l-i}\wedge^{i+j}U\otimes \wedge^{l-(i+j)}W
=\oplus_{s=i}^l \wedge^s U\otimes \wedge^{l-s}W. \]
From this it follows that 
\[ F_i/F_{i+1}=\wedge^i U\otimes \wedge^{l-i}W \]
as $K$-vector space, hence $g$ is an isomorphism of $\lg$-connections,
and the claim of the Lemma follows.
\end{proof}

Let $Op(\K)$ denote the set of all natural transformations
of the underlying set valued functor of $\K$. It follows that $Op(\K)$
is an associative ring. The exterior product $\wedge^l$ from Lemma
\ref{filtration} defines operations $\lambda^l\in Op(\K)$.

\begin{theorem} \label{main2} For any field extension $k\subseteq K$
  of fields of characteristic $p$ there is a contravariant functor
\[ \Kl:Lie_{K/k} \rightarrow \text{\lring} \]
defined by
\[ \Kl(\lg)=(\K(\mod \lg) , \lambda^l_{l\geq 0}) .\]
Here $\K(\mod \lg)$ is the grothendieck ring of the category $\mod
\lg$ and $\lambda^l[V,\rho]=[\wedge^l V,\wedge^l \rho]\in \K(\mod \lg)$.
\end{theorem}
\begin{proof} The proof uses Lemma \ref{filtration} and 
is left to the reader as an exercise.
\end{proof}

\begin{example} Radicial descent. \label{calculation} 
\end{example}
There is an isomorphism of  $\lambda$-rings
\[ \K(\modpf \lg )  \cong \Z. \]
There is by Theorem 2.3 in \cite{cartier} an equivalence between the
  category $\modpf \lg$ of $p$-flat $\lg$-connections of finite dimension and the
  category $\Vect_{K^\lg}$ of $K^\lg$-vectorspaces of finite dimension, hence there is
  an isomorphism of grothendieck groups
\[ e: \K(\modpf \lg )  \cong \K(\Vect_{K^\lg})\cong \Z .\]

In the following we use the operations defined above
to define characteristic classes of arbitrary virtual
connections in the grothendieck ring $\K(\mod \lg )$ as done in
Section 2. Consider the functor
\[ \K:Lie_{K/k}\rightarrow \ring \]
where for any $\lg \in Lie_{K/k}$ $\K(\mod \lg)$ is the grothendieck
ring of the category $\mod \lg$. 

We aim to show that $\K$ has a theory of characteristic classes $c$
using operations defined in the previous section.

Let $u=t/1-t$ and consider the formal powerseries 
\[ \gamma_t(x)=\lambda_u(x)=\sum_{l \geq 0}\lambda^l(x)u^l=\sum_{l\geq
  0}\gamma^l(x)t^l.\]
It has the property that for any exact sequence
\[ 0\rightarrow U\rightarrow V \rightarrow W \rightarrow 0 \]
in $\mod \lg$  there is an equality of formal powerseries in
$1+\K(\mod \lg )[[t]]$:
\[ \gamma_t(V)=\gamma_t(U)\gamma_t(W).\]
For all $l\geq 0$ it follows $\gamma^l \in Op(\K)$.
We get well defined operations -\emph{ gamma operations} -
\[ \gamma^l:\K(\mod \lg )\rightarrow \K(\mod \lg ). \]
satisfying
\[ \gamma^l(V)=\sum_{i+j=l}\gamma^i(U)\gamma^j(W) .\]

\begin{definition} \label{class} Let $x=\sum_in_i[V_i,\rho_i]=\sum_i n_i[V_i]\in
  \K(\mod \lg)$ be an element.
Define the following map
\[ e:\K(\mod \lg)\rightarrow \mathbb{Z} \]
by
\[ e(x)=\sum_i n_i dim_K(V_i).\]
Define furthermore
\[ d:\K(\mod \lg)\rightarrow \K(\mod \lg) \]
by
\[ d(x)=e(x)[\theta_K] \]
where $\theta_K$ is the trivial rank one $\lg$-connection.
For any integer $l\geq 0$ define the \emph{l'th characteristic class}
of $x$ to be
\[ c_l(x)=(-1)^l\gamma^l(x-d(x)) \in \K (\mod \lg).\]
\end{definition}

One checks immediately that the notions above are well defined.

Let $F: \lh \rightarrow \lg$ be a morphism in $Lie_{K/k}$. There is for every $\lg$-connection
$(W,\rho)$ a canonical $\lh$-connection $F^*W=K\otimes W$. This
construction commutes with exterior product: $F^*(\wedge^l
W)=\wedge^l(F^*W)$. We say that an element $x=\sum n_i[V_i] \in \K(\mod \lg)$
is \emph{effective} if $n_i>0$ for all $i$.
We get the following result for the group $\K(\mod \lg )$:

\begin{theorem}\label{mainconn} Let $x,y\in \K(\mod \lg)$ be
  arbitrary elements with $x$ effective.
The following formulas hold:
\begin{align}
c_l(x+y)=\label{class1} &\sum_{i+j=l}c_i(x)c_j(y) \\
c_l(x)=&\label{class2} 0 \text{ if }l> e(x)\\
F^*c_l(x)=\label{class3}&c_l(F^*x)
\end{align}
\end{theorem}

\begin{proof}The proof is similar to the proof of Theorem
  \ref{maintheorem} and is left to the reader as an exercise.
\end{proof}

\begin{definition}\label{total}
Let for any effective $x\in \K (\mod \lg)$
\[c(x)=\sum_{l\geq 0}c_l(x)\in \K(\mod \lg) \]
be it's \emph{virtual characteristic class}.
\end{definition}
By Lemma \ref{mainconn}, Formula \ref{class2} it follows $c(x)$ is
well defined for any $x\in \K (\mod \lg)$.

\section{Differential forms, curves and Ore-extensions}

In this section we make the constructions in Section 2 explicit and
construct characteristic classes for connections on curves and
connections defined in terms of differential forms on fields.

First we introduce the \emph{Cartier operator} (following the presentation
in \cite{cartier} and \cite{katz}) and relate it to the $p$-curvature of a connection
defined on the field $K$. 

Let $K^p\subseteq L \subseteq K$ be fields of
characteristic $p>0$ and let $\lg=\Der_L(K)$ be the $p$-Lie algebra of
derivations of $K$ over $L$. Let $\kp$ be the field of $p$'th roots of
elements of $K$ ie $\kp$ is the splitting field of all polynomials
$T^p-a$ with $a\in K$. It has the property that for all elements $a\in
K$ there is a unique element $x=a^{1/p}\in \kp$ with
$x^p=a$. Furthermore one has that $(a+b)^{1/p}=a^{1/p}+b^{1/p}$. 
Let
$\omega=xdy\in \Omega^1_{K/L}$ be a differential form and define the
following map:
\[ C\omega: \lg \rightarrow \kp \]
by
\[ C\omega (\pa)=(\omega (\pa^p)-\pa^{p-1}(\omega (\pa)) )^{1/p} .\]

\begin{definition} The map $C$ is the Cartier operator for the field
  extension $L\subseteq K$.
\end{definition}

\begin{proposition}The following holds:
\begin{align}
&\label{c1}C(\omega+\omega')=C\omega+ C\omega' \\
&\label{c2}C(x\omega)=x^{1/p}C\omega, x\in L \\
&\label{c3}C(dx)=0
\end{align}
\end{proposition}
\begin{proof} We prove \ref{c1}: 
\[ C(\omega+\omega')(\pa)=(
(\omega+\omega')(\pa)-\pa^{p-1}((\omega+\omega')(\pa)))^{1/p}=\]
\[ (
\omega(\pa^p)-\pa^{p-1}(\omega(\pa))+\omega'(\pa)-\pa^{p-1}(\omega'(\pa))
)^{1/p}= \]
\[
(\omega(\pa)-\pa^{p-1}(\omega(\pa))^{1/p}+(\omega'(\pa)-\pa^{p-1}(\omega'(\pa))
)^{1/p}=C\omega (\pa)+C\omega'(\pa).\]
We prove \ref{c2}: Let $x\in L$. We get
\[ C(x\omega)(\pa)=(x\omega(\pa)-\pa^{p-1}(x\omega(\pa)) )^{1/p}=\]
\[ (x(\omega(\pa)-\pa^{p-1}(\omega (\pa))) )^{1/p}=x^{1/p}C\omega
(\pa).\]
\ref{c3} is obvious. 
\end{proof}

Define the following connection
\[ r:\Der_L(K)\rightarrow \End_L(K) \]
by
\[ r(\pa)(x)=\pa(x)+\omega(\pa)x.\]
Let $R_r(\pa,\pa')=[r(\pa),r(\pa')]-r([\pa,\pa'])$ be the curvature of
$r$ and $\psi_r(\pa)=r(\pa^p)-r(\pa)^p$ the $p$-curvature.

\begin{theorem} \label{cartieroperator} The following holds:
\begin{align}
\label{k1}&R_r(\pa,\pa')=d\omega(\pa,\pa')\\
\label{k2}&\psi_r(\pa)=(-1)^p(C\omega(\pa)-\omega(\pa))^p.
\end{align}
\end{theorem}
\begin{proof} Assume $\omega=xdy$. We first prove \ref{k1}: It follows
\[ d\omega(\pa,\pa')=\pa(x)\pa'(y)-\pa'(x)\pa(y). \]
It follows that 
\[ \pa(\omega(\pa'))-\pa'(\omega(\pa))-\omega([\pa,\pa'])= \]
\[ \pa(x\pa'(y))-\pa'(x\pa(y))-x[\pa,\pa'](y)=\]
\[
\pa(x)\pa'(y)+x\pa\pa'(y)-\pa'(x)\pa(y)-x\pa\pa'(y)-x[\pa,\pa'](y)=\]
\[\pa(x)\pa'(y)-\pa'(x)\pa(y)=d\omega(\pa,\pa').\]
We get thus the formula
\[ \pa(\omega(\pa'))-\pa'(\omega(\pa))=d\omega(\pa,\pa')+\omega([\pa,\pa']).\]
We get
\[
[r\pa,r\pa'](u)=(\pa+\omega(\a))(\pa'+\omega(\pa')(u)-(\pa'+\omega(\pa'))(\pa+\omega(\pa))(u)=\]
\[
\pa\pa(u)+\pa(\omega(\pa')u)+\omega(\pa)\pa'(u)+\omega(\pa)\omega(\pa')u
\]
\[-
\pa'\pa(u)-\pa'(\omega(\pa)u)-\omega(\pa')\pa(u)-\omega(\pa')\omega(\pa)u\]
\[=[\pa,\pa'](u)+\pa(\omega(\pa')u)+\omega(\pa)\pa'(u)-\pa'(\omega(\pa)u)-\omega(\pa')\pa(u)=\]
\[ [\pa,\pa'](u)+\pa(\omega(\pa'))u-\pa'(\omega(\pa))u  =\]
\[ [\pa,\pa'](u)+d\omega(\pa,\pa')(u)+\omega([\pa,\pa'])(u). \]
It follows that
\[ R_r(\pa,\pa')=[r\pa,r\pa']-r([\pa,\pa'])=d\omega(\pa,\pa')\]
hence \ref{k1} is proved.
In the ring $\End_L(K)$ the following holds:
\[ (a+\pa)^p=a^p+\pa^p+\pa^{p-1}(a). \]
This is proved using induction. 
We prove \ref{k2}: 
\[
\psi_r(\pa)=r(\pa^p)-r(\pa)^p=\pa^p+\omega(\pa^p)-(\pa+\omega(\pa))^p=\]
\[ \pa^p+\omega(\pa^p)-\omega(\pa)^p-\pa^p-\pa^{p-1}(\omega(\pa)) =\]
\[ \omega(\pa^p)-\pa^{p-1}(\omega(\pa))-\omega(\pa)^p= \]
\[ C\omega(\pa)^p-\omega(\pa)^p=(C\omega(\pa)-\omega(\pa))^p,\]
and the claim follows.

\end{proof}

\begin{example} Differential forms.  \label{mainexample} 
\end{example}
Let $\w\in \Omega=\Omega^1_{K/L}$ and let $\lg=\Der_L(K)$
Define the following map
\[ \rho_\w: \lg \rightarrow \End_L(K) \]
by 
\[ \rho_\w(\pa)(x)=\pa(x)+\w(\pa)(x) .\]
It follows $\rho_\w$ is a connection and from Proposition
\ref{cartieroperator} one gets
\[ R_{\rho_\w}(\pa,\pa')=d\w(\pa,\pa') \]
and
\[ \psi_{\rho_\w}(\pa)=(-1)^p(C\w(\pa)-\w(\pa))^p ,\]
where $C\w$ is the Cartier operator (see \cite{cartier} and \cite{katz}).
Note: By Proposition 7 in \cite{cartier} it follows 
that $R_{\rho_\w}=\psi_{\rho_\w}=0$ if and only if
$\w=dlog(x)=x^{-1}dx$.
We get thus a map
\[ \phi:\Omega \rightarrow \K (\mod \lg ) \]
defined by
\[ \phi(\w)=c(K,\rho_\w) \]
detecting if a differential form $\w$ is a logarithmic derivative.

\begin{example} \label{adjoint} Adjoint representation.
\end{example}
 Let
  $\lg\subseteq \Der_k(K)$ be a restricted Lie-Rinehart algebra. There is a representation
\[ ad: \lg \rightarrow \End_k(\Der_k(K)) \]
defined by
\[ ad(x)(y)=[x,y] .\]
One sees
\[ ad(x)(\alpha y)=[x,\alpha y]=x(\alpha y)-\alpha yx=x(\alpha
)y+\alpha xy-\alpha yx=\alpha [x,y]+x(\alpha )y =\]
\[ \alpha ad(x)(y)+x(\alpha)y ,\]
hence $ad$ is a connection. 
The Jacobi-identity shows $R_{ad}=0$ hence the map $ad$ is a representation of Lie
algebras. 
In $\End_k(\Der_k(K))$ the following formula holds:
\[ [\pa[\pa[\cdots [\pa,\eta]\cdots
]=\sum_{i=0}^n\binom{n}{i}\pa^{n-i}\eta \pa^i\]
for all $n\geq 1$.
We get the formula
\[
ad(\pa^p)(\eta)=[\pa^p,\eta]=\pa^p\eta-\eta\pa^p=\sum_{i=0}^p\binom{p}{i}\pa^{p-i}\eta
\pa^i =\]
\[ [\pa[\pa[\cdots [\pa,\eta]\cdots ]=ad(\pa)^p(\eta) \]
hence it follows
\[ \psi_\rho(\pa)=ad(\pa^p)-ad(\pa)^p=0\]
and it follows $ad$ is a flat $\lg$-connection on $\Der_k(K)$
with zero $p$-curvature.

One may also check that for any ideal $I\subseteq \Der_k(K)$ it
follows $I$ is closed under $p$-powers, hence $I$ is a restricted Lie-Rinehart algebra.
The map
\[ad:\Der_k(K)\rightarrow \End_k(I) \]
and
\[ ad:\Der_k(K)\rightarrow\End_k(\Der_k(K)/I) \]
makes $I$ and $\Der_k(K)/I$ into $p$-representations and the sequence
\[ 0\rightarrow I \rightarrow \Der_k(K) \rightarrow
\Der_k(K)/I\rightarrow 0 \]
is an exact sequence of $p$-flat $\lg$-connections.
We get a characteristic class
\[ c(\Der_k(K)/I, ad)\in \K(\mod \lg ) \]
for each ideal $I\subseteq \Der_k(K)$.

There is a $p$-flat connection 
\[ ad:\lg \rightarrow \End_k(U^{[p]}(\lg) ) \]
where $U^{[p]}(\lg)$ is the restricted enveloping algebra of $\lg$, hence we get a characteristic class 
\[ c(U^{[p]}(\lg),ad)\in \K(\mod \lg ).\]

\begin{example}Curves and Ore extensions.
\end{example}

If
$\pi:C\rightarrow C'$ is a finite morphism of projective curves
over an arbitrary field $k$ and $K'\subseteq K$ is 
the corresponding finite extension of function fields, we get an exact 
sequence of Lie-Rinehart algebras

\[ 0\rightarrow \lh \rightarrow \Der_k(K)\rightarrow^{d\pi}
K\otimes_{K'}\Der_k(K')   \rightarrow 0.\]
The Lie-Rinehart algebra $\lh=\Der_{K'}(K)$ is a finite dimensional 
$K$-vector space. Let $\lg =\Der_k(K)$ and $\ll=\Der_k(K')$. 
The Lie-Rinehart algebras $\lg$ and $\ll$ are infinite dimensional
$k$-Lie algebras.
We get for any flat $\lg$-connection $(V,\rho)$
and $i\geq 0$ canonical flat connections - \emph{Gauss-Manin
  connections} (see \cite{maakestad1}) -
\[ \nabla_{GM}^{i,\rho}:\ll \rightarrow \End_k(\H^i(\lh, V)) \]
where $\H^i(\lh, V)$ is the Chevalley-Hochschild cohomology of $V$ as
$\lh$-module. Let $\K (\modf \ll)$ denote the grothendieck ring of the
category $\modf \ll$.
We get a well defined cohomology class
\[ d\pi_!(V,\rho)=\sum_{i\geq 0}(-1)^i[\H^i(\lh,V)]\in \K(\modf \ll ) .\]
The connection $d\pi^i(\rho)$ is not $p$-flat in general.
(See \cite{katz0} for examples where the $p$-curvature of
$d\pi^i(\rho)$ is related to the Kodaira-Spencer class).

Let $\rho=ad$ from Example \ref{adjoint} with 
\[ \rho: \Der_k(K)\rightarrow \End_k(\Der_k(K)) \]
and make the following definition:
\begin{definition}
Let the cohomology class
\[ \Delta (\pi)=d\pi_!(\Der_k(K), \rho)=\sum_{i\geq 0}(-1)^i[\H^i(\lh, \Der_k(K))]\in
\K(\modf \ll ) \]
be the \emph{ramification class} of the morphism $\pi$.
\end{definition}
If $\pi$ is a separable morphism, it follows $\lh=0$. In this case
\[\Delta (\pi) =[\Der_k(K)]\in \K(\modf \ll ).\]
Hence the cohomology class 
\[ \Delta (\pi) \in \K(\modf \ll ) \]
is related to the ramification of the morphism $\pi$. 
We get a characteristic class
\[ c(\pi)=\sum_{i\geq 0}c_i(\Delta(\pi))\in \K (\modf \ll)  \]
defined for an arbitrary finite morphism $\pi:C\rightarrow C'$  of
projective curves.

 Assume there is a separable morphism of 
curves
\[ C' \rightarrow \P^1_k \]
where $k$ is a field of characteristic zero and let $(V,\rho)$ be a
flat $\lg$-connection.
The function field $K(\P^1_k)$ equals $k(t)$ where $t$ is a transcendental
variable over $k$. We get an isomorphism
\[ \ll\cong K'\pa \]
where $\pa$ is partial derivative with respect to $t$. It follows
that the category of flat finite dimensional $\ll$-connections
is equivalent to the category of $K'\{T\}$-modules of finite dimension
over $K'$. Here $K'\{T\}$ is the \emph{Ore extension} of $K'$ by $T$. It is the
twisted polynomial ring on the variable $T$ with multiplication defined as follows:
\[ Ta=aT+\pa(a) \]
for any $a\in K'$. In this case the cohomology group $\H^i(\lh, V)$
is a finite dimensional $K$-vector space. It follows $\H^i(\lh,V)$ is
a finite dimensional $K'$-vector space. The Gauss-Manin connection
\[ \nabla^{i,\rho}_{GM}: \ll\rightarrow \End_k(\H^i(\lh,V)) \]
makes $\H^i(\lh,V)$ into a finite dimensional left $K'\{T\}$-module. 
By the cyclic vector theorem (see \cite{kedlaya} Theorem 5.4.2) 
any finite dimensional left $K'\{T\}$-module is on the form
\[ K'\{T\}/P(T)K'\{T\} \]
where $P(T)\in K'\{T\}$ is a non-commutative polynomial.
It follows there is a non-commutative polynomial
\[ P_i(T)\in K'\{T\} \]
such that
\[ \H^i(\lh,V)\cong K'\{T\}/P_i(T)K'\{T\} \]
as $K'\{T\}$-module.
In this case we get a cohomology class
\[ \Delta(V)=\sum_{i\geq 0}(-1)^i[\H^i(\lh,V)]=\sum_{i\geq
  0}(-1)^i[K'\{T\}/P_i(T)K'\{T\}]\in \K(\modf \ll).\]
We get a characteristic class
\[ c(V)=\sum_{i\geq 0}c_i(\Delta(V))\in \K(\modf \ll) .\]
By a Theorem of Quillen (see Theorem 1 in \cite{mcconnell}) it follows
\[ \hK_n(\modf \ll)\cong \hK_n(K'\{T\})\cong \hK_n(K')\cong \mathbb{Z} \]
for all $n\geq 0$. It follows the class
\[ [\H^i(\lh,V),\nabla_{GM}^{i,\rho}]\in \K(\modf \ll) \]
is independent of choice of connection whenever $C'$ is a separable
cover of $\P^1_k$. It follows that the grothendieck group
$\K(\modf \ll)$ does not contain enough information to detect invariants
of the connection $\nabla_{GM}^{i,\rho}$.

\textbf{Acknowledgements}. Thanks to Ken Goodearl, Max Karoubi and Charles Weibel for discussions and
comments.

\end{document}